\documentclass[12pt]{article}
\usepackage{amsmath,amssymb,amsbsy,amscd}
\newtheorem{thm}{Theorem}[section]
\newtheorem{prop}[thm]{Proposition}
\newtheorem{lem}[thm]{Lemma}
\newtheorem{cor}[thm]{Corollary}

\newtheorem{pf}{\it Proof.}

\newcommand{\qed}{{\hfill $\Box $}}          
\newcommand{\Aut }{\mathop{\rm Aut}\nolimits}
\newcommand{\lcm}{\mathop{\rm lcm}\nolimits}
\newcommand{\zs}{\{ 0\} }
\newcommand{\sm}{\setminus}

\newcommand{\N}{{\bf N}}
\newcommand{\Z}{{\bf Z}}
\newcommand{\Zn}{{\Z _{\geq 0}}}
\newcommand{\x}{{\bf x}}

\newcommand{\w}{{\bf w}}
\newcommand{\kx}{k[{\bf x}]}

\begin{document}

\title{A generalization of the Shestakov-Umirbaev inequality} 

\author{Shigeru Kuroda}

\date{}

\maketitle

\section{Introduction}\setcounter{equation}{0}

Let $k$ be a field, and 
$\kx =k[x_1,\ldots ,x_n]$ the polynomial ring in $n$ variables 
over $k$ for $n\in \N $. 
For a polynomial $\Phi =\sum _{i=0}^l\phi _iy^i$ 
in a variable $y$ over $\kx $ and $g\in \kx $, 
we denote $\Phi (g)=\sum _{i=0}^l\phi _ig^i$, 
where $\phi _0,\ldots ,\phi _l\in \kx $ for $l\geq 0$. 
Then, it follows that 
$$
\deg ^g\Phi :=\max \{ \deg (\phi _ig^i)\mid i=0,\ldots ,l\} \geq 
\deg \Phi (g)
$$
in general. 
Here, 
$\deg f$ denotes the total degree of $f$ 
for each $f\in \kx $. 
Shestakov-Umirbaev~\cite[Theorem 3]{SU1} proved 
an inequality which describes the difference 
between $\deg ^g\Phi $ and $\deg \Phi (g)$. 
Using this result, 
they settled in~\cite{SU2} an important open problem on 
automorphisms of $\kx $ as follows.

Let $\sigma :\kx \to \kx $ be a homomorphism of $k$-algebras. 
Then, 
$\sigma $ is an isomorphism if and only if 
\begin{equation}\label{eq:isom}
k[\sigma (x_1),\ldots ,\sigma (x_n)]=\kx . 
\end{equation}
For example, 
$\sigma $ is an isomorphism if there exist 
$(a_{i,j})_{i,j}\in GL_n(k)$ 
and $(b_i)_i\in k^n$ such that 
$\sigma (x_i)=\sum _{j=1}^na_{i,j}x_j+b_i$ for each $i$. 
It also follows that 
$\sigma $ is an isomorphism if there exists 
$l\in \{ 1,\ldots ,n\} $ such that 
$\sigma (x_i)=x_i$ for each $i\neq l$ 
and $\sigma (x_l)=\alpha x_l+f$ 
for some $\alpha \in k^{\times }$ 
and $f\in k[x_1,\ldots ,x_{l-1},x_{l+1},\ldots ,x_n]$. 
An automorphism of $\kx $ as in the former example 
is said to be {\it affine}, 
and one as in the latter example 
is said to be {\it elementary}. 
Because an invertible matrix 
is expressed as a product of elementary matrices, 
each affine automorphism 
can be obtained by the composition of elementary automorphisms. 
Then, 
a problem arises whether the automorphism group $\Aut _k\kx $ 
can be generated by elementary automorphisms. 
This is called the {\it Tame Generators Problem}. 
If $n=1$, then every automorphism of $\kx $ is in fact elementary. 
If $n=2$, 
then $\Aut _k\kx $ is generated by elementary automorphisms, 
which was shown by 
Jung~\cite{Jung} in 1942 in case $k$ is of characteristic zero, 
and by van der Kulk~\cite{Kulk} in 1953 for an arbitrary $k$. 
We note that this result is a consequence 
of the following characterization of automorphisms of $\kx $. 

\begin{prop}\label{prop:Jung}
If $n=2$, then 
either $\deg \sigma (x_1)|\deg \sigma (x_2)$ or 
$\deg \sigma (x_2)|\deg \sigma (x_1)$ holds 
for each $\sigma \in \Aut _k\kx $. 
\end{prop}
Here, $a|b$ denotes that $b$ is divisible by $a$ 
for each $a,b\in \N $. 
Due to (\ref{eq:isom}), 
$\deg \sigma (x_i)$ must be positive for each 
$\sigma \in \Aut _k\kx $ and $i=1,\ldots ,n$.

When $n\geq 3$, 
the problem becomes extremely difficult. 
In 1972, 
Nagata~\cite{Nagatame} conjectured that 
the automorphism $\tau \in \Aut _k\kx $ for $n=3$ defined by 
\begin{equation*}\label{eq:nagata}
\tau (x_1)=x_1-2(x_1x_3+x_2^2)x_2-(x_1x_3+x_2^2)^2x_3,\  
\tau (x_2)=x_2+(x_1x_3+x_2^2)x_3,\ 
\tau (x_3)=x_3
\end{equation*}
cannot be obtained by the composition of 
elementary automorphisms of $\kx $. 
This conjecture was well-known, 
but was not settled for a long time. 
In 2004, however, 
Shestakov-Umirbaev~\cite{SU2} 
finally showed that the Nagata conjecture is true 
if $k$ is of characteristic zero. 
The inequality mentioned at the beginning 
plays a crucial role 
in their solution of the Nagata conjecture. 
The Tame Generators Problem is thus settled for $n=3$, 
but remains open for $n\geq 4$. 
We note that 
the extension $\tilde{\tau }\in \Aut _k\kx $ 
of the Nagata automorphism $\tau $ for $n\geq 4$ 
defined by $\tilde{\tau }(x_i)=\tau (x_i)$ for $i=1,2,3$ 
and $\tilde{\tau }(x_i)=x_i$ for $i=4,\ldots ,n$ 
is a composite of elementary automorphisms (see~\cite{Smith}).

The argument in~\cite{SU2} is indeed difficult, 
but employs no advanced facts other than those in~\cite{SU1}. 
Therefore, 
the results in~\cite{SU1} are of great importance. 
However, 
its argument is also difficult, 
and, consequently, the proof of 
this landmark work of Shestakov-Umirbaev 
is unfortunately not widely understood.

The purpose of the present paper 
is to generalize the results of~\cite{SU1}. 
Our argument is quite simple and elementary, 
but the results are general and very interesting. 
These results will be useful 
not only for a better understanding 
of the theory of Shestakov-Umirbaev, 
but also to generalize it to higher dimensions 
to solve the Tame Generators Problem for $n\geq 4$. 
As an application, 
we give a generalization of Proposition~\ref{prop:Jung} 
in Theorem~\ref{thm:general-Jung}.

Section~\ref{sect:basic} is devoted to proving a basic result. 
We derive its consequence in Section~\ref{sect:SUineq}, 
and apply it to characterizations of automorphisms of $\kx $ 
in Section~\ref{sect:auto}. 
In Section 5, 
we generalize a lemma~\cite[Lemma 5]{SU1} 
of Shestakov-Umirbaev which also plays an important role 
in the solution of the Nagata conjecture.

It should be noted that 
Makar-Limanov~\cite{EMW} also gave another proof of~\cite[Theorem 3]{SU1} 
in a different fashion.

\section{Differentials}
\setcounter{equation}{0}\label{sect:basic}

In what follows, 
we always assume that $k$ is of characteristic zero. 
First, 
we introduce some terminology concerning the grading of a polynomial ring.

Let $\Gamma $ be a totally ordered additive group, 
and $\w =(w_1,\ldots ,w_n)$ an element of $\Gamma ^n$. 
We define the $\w $-{\it weighted grading} 
$\kx =\bigoplus _{\gamma \in \Gamma }\kx _{\gamma }$  
by setting $\kx _{\gamma }$ 
to be the $k$-vector space generated by 
$x_1^{a_1}\cdots x_n^{a_n}$ for $a_1,\ldots ,a_n\in \Zn $ 
with $\sum _{i=1}^na_iw_i=\gamma $ for each $\gamma \in \Gamma $. 
Here, $\Zn $ denotes the set of nonnegative integers, 
and $l\gamma $ denotes the sum of $l$ copies of $\gamma $ 
for each $l\in \Zn $ and $\gamma \in \Gamma $. 
It follows that 
$\kx _{\gamma }\kx _{\mu }\subset \kx _{\gamma +\mu }$ 
for each $\gamma ,\mu \in \Gamma $. 
Assume that $f=\sum _{\gamma \in \Gamma }f_{\gamma }$ 
is an element of $\kx $, 
where $f_{\gamma }\in \kx _{\gamma }$ for each $\gamma $. 
If $f\neq 0$, 
then the $\w $-{\it degree} $\deg _{\bf w}f$ of $f$ 
is defined to be the maximum among 
$\gamma \in \Gamma $ with $f_{\gamma }\neq 0$. 
If $f=0$, then we set $\deg _{\bf w}f=-\infty $, 
i.e., 
a symbol which is less than each element of $\Gamma $. 
The addition is defined by 
$(-\infty )+\gamma =\gamma +(-\infty )=-\infty $ 
for each $\gamma \in \Gamma \cup \{ -\infty \} $, 
and the sum of $l$ copies of $-\infty $ is 
denoted by $l(-\infty )$ for each $l\in \Zn $. 
We say that $f$ is $\w $-{\it homogeneous} 
if $f=f_{\gamma }$ for some $\gamma $. 
In case $f\neq 0$, 
we define $f^{\w }=f_{\delta }$, 
where $\delta =\deg _{\w }f$. 
Then, it follows that 
$\deg _{\w }f^{\w }=\deg _{\w }f$, 
$\deg _{\w }(f-f^{\w })<\deg _{\w }f$, 
and $(f_1f_2)^{\w }=f_1^{\w }f_2^{\w }$ 
for each $f,f_1,f_2\in \kx \sm \zs $. 
We denote by $\Gamma _{\geq 0}$ 
the set of $\gamma \in \Gamma $ 
with $\gamma \geq 0$, 
where $0$ is the zero of the additive group $\Gamma $. 
We remark that $\deg _{\w }f\geq 0$ holds 
for each $f\in \kx \sm \zs $ 
whenever $\w $ is an element of $(\Gamma _{\geq 0})^n$. 
If $\Gamma =\Z $ and $\w =(1,\ldots ,1)$, 
then the $\w $-degree is the same as the total degree.

Now, 
for $\Phi \in \kx [y] $ and $g\in \kx $, 
we define 
\begin{equation}\label{eq:deg_w^g}
\deg _{\w }^g\Phi =\max \{ \deg _{\w }(\phi _ig^i)\mid i\in \Zn \} , 
\end{equation}
where $\phi _i\in \kx $ for each $i\in \Zn $ with $\Phi =\sum _i\phi _iy^i$. 
Then, 
$\deg _{\w }^g\Phi $ is at least $\deg _{\w }\Phi (g)$ in general. 
The purpose of this section is to prove an inequality which 
describes the difference 
between $\deg _{\w }\Phi (g)$ and $\deg _{\w }^g\Phi $.

Let $\partial _y^i\Phi $ denote the $i$-th order derivative of $\Phi $ 
in $y$ for each $i\in \Zn $, 
and $\deg _y\Phi $ the degree of $\Phi $ in $y$. 
Obviously, 
$\deg _{\w }^g\partial _y^i\Phi 
=\deg _{\w }\!\left(\partial _y^i\Phi \right)\!(g)$ 
if $i\geq \deg _y\Phi $. 
So, 
we may consider the nonnegative integer 
\begin{equation}
m_{\w }^g(\Phi )=\min \left\{ i\in \Zn \mid
\deg _{\w }^g\partial _y^i\Phi 
=\deg _{\w }\!\left(\partial _y^i\Phi \right)\!(g)\right\} . 
\end{equation}
Note that 
\begin{equation}\label{eq:ref3}
m_{\w }^g(\Phi )=
m_{\w }^g(\partial _y\Phi )+1\quad 
\text{ and } \quad 
\deg _{\w }^g\Phi =
\deg _{\w }^g\partial _y\Phi +\deg _{\w }g 
\end{equation}
if $m_{\w }^g(\Phi )\geq 1$ and $g\neq 0$, 
where $\partial _y\Phi =\partial _y^1\Phi $, 
since $k$ is of characteristic zero.

Let $\Omega _{\kx /k}$ be the differential module of 
$\kx $ over $k$, 
and $\bigwedge ^r\Omega _{\kx /k}$ 
the $r$-th exterior power of the 
$\kx $-module $\Omega _{\kx /k}$ 
for $r\in \{ 1,\ldots ,n\} $. 
Then, 
each $\omega \in \bigwedge ^r\Omega _{\kx /k}$ 
is uniquely expressed as 
$$
\omega =\sum _{1\leq i_1<\cdots <i_r\leq n}
f_{i_1,\ldots ,i_r}dx_{i_1}\wedge \cdots \wedge dx_{i_r},
$$
where $f_{i_1,\ldots ,i_r}\in \kx $ for each $i_1,\ldots ,i_r$. 
Here, 
$df$ denotes the differential of $f$ 
for each $f\in \kx $. 
We define the $\w $-degree of $\omega $ by 
\begin{equation}\label{eq:def:deg}
\deg _{\w }\omega =\max \{ \deg _{\w }(f_{i_1,\ldots ,i_r})
+w_{i_1}+\cdots +w_{i_r}
\mid 1\leq i_1<\cdots <i_r\leq n\} . 
\end{equation}
Since $df=\sum _{i=1}^n(\partial f/\partial x_i)dx_i$ 
and $k$ is of characteristic zero, 
the equality 
\begin{equation}\label{eq:ref0}
\deg _{\w }df=\max \left\{ \deg _{\w }\!\left(
\frac{\partial f}{\partial x_i}\right) +w_i
\mid i=1,\ldots ,n\right\} =\deg _{\w }f
\end{equation}
holds for each $f\in \kx $. 
It is easily verified that 
$\deg _{\w }(\omega +\omega ')\leq 
\max \{ \deg _{\w }\omega ,\deg _{\w }\omega '\} $, 
\begin{equation}\label{eq:ref1}
\deg _{\w }(\omega \wedge \eta )\leq \deg _{\w }\omega +\deg _{\w }\eta  \ \ \ \
\text{and} \ \ \ 
\deg _{\w }(f\omega )=\deg _{\w }f+\deg _{\w }\omega 
\end{equation}
for each 
$\omega ,\omega '\in \bigwedge ^r\Omega _{\kx /k}$ 
and $\eta \in \bigwedge ^s\Omega _{\kx /k}$ 
for $r,s\in \{ 1,\dots ,n\} $ with $r+s\leq n$, 
and  $f\in \kx $. 

In the notation above, we have the following theorem.

\begin{thm}\label{thm:main}
Let $f_1,\ldots ,f_r$ be elements of $\kx $ for $r\geq 1$ 
which are algebraically independent over $k$, 
and set $\omega =df_1\wedge \cdots \wedge df_r$. 
Then, the inequality 
\begin{equation}\label{eq:ref5}
\deg _{\w }\Phi (g)\geq \deg _{\w }^g\Phi +m_{\w }^g(\Phi )
(\deg _{\w }(\omega \wedge dg)-\deg _{\w }\omega -\deg _{\w }g) 
\end{equation}
holds for each $\Phi \in k[f_1,\ldots ,f_r][y]\sm \zs $, 
$g\in \kx \sm \zs $ and $\w \in \Gamma ^n$. 
\end{thm}
\begin{pf}\rm
Recall that, 
for $h_1,\ldots ,h_s\in \kx $ for $s\geq 1$, 
it follows that $h_1,\ldots ,h_s$ 
are algebraically independent over $k$ 
if and only if $dh_1\wedge \cdots \wedge dh_s\neq 0$ 
when $k$ is of characteristic zero. 
Therefore, 
$\omega \wedge df_i=0$ for $i=1,\ldots ,r$, 
while $\omega \neq 0$ by assumption. 
By chain rule, we may write 
$d(\Phi (g))=(\partial _y\Phi )(g)dg+\sum _{i=1}^r\psi _idf_i$, 
where $\psi _i\in \kx $ for each $i$. 
Thus, 
\begin{align}\begin{split}\label{eq:proofref1}
\omega \wedge d(\Phi (g))=
(\partial _y\Phi )(g)\omega \wedge dg
+\sum _{i=1}^r\psi _i\omega \wedge df_i
=(\partial _y\Phi )(g)\omega \wedge dg. 
\end{split}\end{align}
By (\ref{eq:ref0}), (\ref{eq:ref1}) and (\ref{eq:proofref1}), we have
\begin{align}\begin{split}\label{eq:proofref2}
&\deg _{\w }\omega +\deg _{\w }\Phi (g)
=\deg _{\w }\omega +\deg _{\w }d(\Phi (g))
\geq \deg _{\w }(\omega \wedge d(\Phi (g)))\\
&\qquad \qquad \qquad \quad 
=\deg _{\w }((\partial _y\Phi )(g)\omega \wedge dg) 
=\deg _{\w }(\partial _y\Phi )(g) +
\deg _{\w }(\omega \wedge dg). 
\end{split}\end{align}
By adding $-\deg _{\w }\omega $ 
to both sides of (\ref{eq:proofref2}), 
we get 
\begin{equation}\label{eq:ref6}
\deg _{\w }\Phi (g)\geq 
\deg _{\w }(\partial _y\Phi )(g)+
\deg _{\w }(\omega \wedge dg) -\deg _{\w }\omega . 
\end{equation}
Now, we show (\ref{eq:ref5}) by induction on $m_{\w }^g(\Phi )$. 
If $m_{\w }^g(\Phi )=0$, 
then $\deg _{\w }\Phi (g)=\deg _{\w }^g\Phi $ 
by the definition of $m_{\w }^g(\Phi )$. 
In this case, 
(\ref{eq:ref5}) is clear. 
Assume that $m_{\w }^g(\Phi )\geq 1$. 
Then, 
$m_{\w }^g(\partial _y\Phi )$ is less than $m_{\w }^g(\Phi )$ 
by (\ref{eq:ref3}). 
By induction assumption, and by the equalities in (\ref{eq:ref3}), 
we obtain 
\begin{equation}\label{eq:inductionstep2}
\deg _{\w }(\partial _y\Phi )(g) \geq \deg _{\w }^g\partial _y\Phi 
+m_{\w }^g(\partial _y\Phi )M
=(\deg _{\w }^g\Phi -\deg _{\w }g)+(m_{\w }^g(\Phi )-1)M,
\end{equation}
where $M=\deg _{\w }(\omega \wedge dg)-\deg _{\w }\omega -\deg _{\w }g$. 
Using (\ref{eq:ref6}) and (\ref{eq:inductionstep2}), 
we arrive at 
\begin{align*}
\deg _{\w }\Phi (g)&\geq 
\deg _{\w }(\partial _y\Phi )(g) +
\deg _{\w }(\omega \wedge dg) -\deg _{\w }\omega \\ 
&\geq (\deg _{\w }^g\Phi -\deg _{\w }g)+(m_{\w }^g(\Phi )-1)M
+\deg _{\w }(\omega \wedge dg) -\deg _{\w }\omega \\
&=\deg _{\w }^g\Phi +m_{\w }^g(\Phi )
(\deg _{\w }(\omega \wedge dg)-\deg _{\w }\omega -\deg _{\w }g). 
\end{align*}
Therefore, the inequality (\ref{eq:ref5}) is true. 
\qed\end{pf}

\section{The Shestakov-Umirbaev inequality}
\setcounter{equation}{0}
\label{sect:SUineq}
In this section, we derive some consequences of Theorem~\ref{thm:main}.

First, we remark that 
the element $\deg _{\w }^g\Phi $ of $\Gamma $ 
defined as in (\ref{eq:deg_w^g}) 
is equal to the $(\w ,\deg _{\w }g)$-degree of $\Phi $ 
for each $\Phi \in \kx [y]\sm \zs $, $g\in \kx \sm \zs $ 
and $\w \in \Gamma $, 
where we regard $\Phi $ 
as a polynomial in the $n+1$ variables $x_1,\ldots ,x_n$ and $y$ 
over $k$. 
We denote $\Phi ^{(\w ,\deg _{\w }g)}$ by $\Phi ^{\w ,g}$, 
for short.

\begin{lem}\label{lem:initial-equiv}
Let 
$\Phi \in \kx [y]\sm \zs $, 
$g\in \kx \sm \zs $ and $\w \in \Gamma $. 

{\rm (i)} 
The following conditions are equivalent: 

\quad {\rm (1)} $m_{\w }^g(\Phi )=0$. 

\quad {\rm (2)} $\deg _{\w }^g\Phi =\deg _{\w }\Phi (g)$. 

\quad {\rm (3)} $\Phi ^{\w ,g}(g^{\w })\neq 0$. 

\quad {\rm (4)} $\Phi (g)\neq 0$ and $\Phi (g)^{\w }=\Phi ^{\w ,g}(g^{\w })$. 

{\rm (ii)}  It follows that 
$m_{\w }^g(\Phi )=\min \left\{ i\in \Zn \mid 
\left(\partial _y^i(\Phi ^{\w ,g})\right)\!(g^{\w })\neq 0\right\} $.
\end{lem}
\begin{pf}\rm
(i) 
The equivalence between (1) and (2) 
immediately follows from 
the definition of $m_{\w }^g(\Phi )$. 
In the following, 
we will establish that 
\begin{equation}\label{eq:pequiv1}
\deg _{\w }(\Phi (g)-\Phi ^{\w ,g}(g^{\w }))<\deg _{\w }^g\Phi. 
\end{equation}
Assuming this, 
we can readily check that (2), (3) and (4) are equivalent, since 
$$
\Phi (g)=\Phi ^{\w ,g}(g^{\w })+(\Phi (g)-\Phi ^{\w ,g}(g^{\w })), 
$$ 
and $\Phi ^{\w ,g}(g^{\w })$ 
is contained in $\kx _{\delta }$, where $\delta =\deg _{\w }^g\Phi $.

Write $\Phi =\sum _{i}\phi _iy^i$ 
and $\Phi ^{\w ,g}=\sum _{i}\phi _i'y^i$, 
where $\phi _i,\phi _i'\in \kx $ for each $i$. 
Then, 
$\deg _{\w }(\phi _ig^i)\leq \deg _{\w }^g\Phi $ for each $i$. 
Note that 
$\phi _i'=\phi _i^{\w }$ 
if $\deg _{\w }(\phi _ig^i)=\deg _{\w }^g\Phi$, 
and $\phi _i'=0$ otherwise. 
We have 
$$\phi _ig^i-\phi _i'(g^{\w })^i
=\phi _ig^i-\phi _i^{\w }(g^{\w })^i
=\phi _ig^i-(\phi _ig^i)^{\w }$$ 
in the former case, 
and 
$\phi _ig^i-\phi _i'(g^{\w })^i=\phi _ig^i$ 
in the latter case. 
In each case, 
$\deg _{\w }^g\Phi $ is greater than 
the $\w $-degree of 
$\phi _ig^i-\phi _i'(g^{\w })^i$, 
and hence 
greater than that of 
$$
\sum _{i}\left(\phi _ig^i-\phi _i'(g^{\w })^i\right)
=\Phi (g)-\Phi ^{\w ,g}(g^{\w }).
$$
Thus, we obtain (\ref{eq:pequiv1}), 
thereby proving that (2), (3) and (4) are equivalent.

(ii) Observe that 
$\left(\partial _y^i\Phi \right)^{\w ,g}=\partial _y^i(\Phi ^{\w ,g})$ 
for each $i\in \Zn $. 
In view of this equality, 
it follows that 
$\deg _{\w }^g \partial _y^i\Phi 
=\deg _{\w }\!\left(\partial _y^i\Phi \right)\!(g)$ 
if and only if 
$\left(\partial _y^i(\Phi ^{\w ,g})\right)\!(g^{\w })\neq 0$ 
by the equivalence between (2) and (3) in (i). 
Then, 
the assertion immediately follows from 
the definition of $m_{\w }^g(\Phi )$. 
\qed\end{pf}

Now, 
let $A$ be a $k$-subalgebra of $\kx $, 
and $K$ the field of fractions of $A$. 
We define the {\it initial algebra} $A^{\w }$ of $A$ for $\w $ 
to be the $k$-subalgebra of $\kx $ 
generated by $f^{\w }$ for $f\in A\sm \zs $. 
Then, 
$\Phi ^{\w ,g}$ belongs to $A^{\w }[y]\sm \zs $ 
for each $\Phi \in A[y]\sm \zs $ for any $g\in \kx \sm \zs $. 
We claim that the field of fractions of $B^{\w }$ 
is equal to that of $A^{\w }$ 
whenever $B$ is a $k$-subalgebra of $\kx $ 
whose field of fractions is equal to $K$. 
Indeed, 
if $fg_1=g_2$ for $f\in A$ (resp.\ $f\in B$) 
and $g_1,g_2\in B$ (resp.\ $g_1,g_2\in A$), 
then we have $f^{\w }g_1^{\w }=(fg_1)^{\w }=g_2^{\w }$, 
so $f^{\w }$ belongs to the field of fractions of 
$B^{\w }$ (resp.\ $A^{\w }$). 
For this reason, 
we denote 
the field of fractions of $A^{\w }$ by $K^{\w }$.

For an integral domain $R$ and an element $s$ 
of an integral domain $S$ containing $R$, 
we define $I(R,s)$ to be 
the kernel of the substitution map 
$R[y]\ni f\mapsto f(s)\in S$. 
When $I(R,s)$ is a principal ideal of $R[y]$, 
a generator of $I(R,s)$, 
which is unique up to multiplication by units in $R$, 
is denoted by $P(R,s)$. 
We remark that 
$I(R,s)$ is always principal if $R$ is a unique factorization domain. 
If $R$ is a field and $s$ is algebraic over $R$, 
then we may take $P(R,s)$ to be the minimal polynomial of $s$ over $R$.

\begin{prop}\label{prop:m_g}
Let $A$ be a $k$-subalgebra of $\kx $, 
and $K$ the field of fractions of $A$. 
Then, for each 
$\Phi \in A[y]\sm \zs $, $g\in \kx \sm \zs $ and $\w \in \Gamma ^n$, 
we have the following:

{\rm (i)} If $g^{\w }$ is transcendental over $K^{\w }$, 
then $m_{\w }^g(\Phi )=0$ and 
$\deg _{\w }\Phi (g)=\deg _{\w }^g\Phi $. 

{\rm (ii)} 
If $g^{\w }$ is algebraic over $K^{\w }$, 
then $m_{\w }^g(\Phi )$ is at most the quotient of 
$\deg _y\Phi ^{\w ,g}$ divided by 
$[K^{\w }(g^{\w }):K^{\w }]$. 
If furthermore $I(A^{\w },g^{\w })$ is a principal ideal, 
then there exists $H\in A^{\w }[y]\sm I(A^{\w },g^{\w })$ 
such that 
$\Phi ^{\w ,g}=P(A^{\w },g^{\w })^{m}H$, 
where $m=m_{\w }^g(\Phi )$. 
\end{prop}
\begin{pf}\rm
(i) If $g^{\w }$ is transcendental over $K^{\w }$, 
then $\Phi ^{\w ,g}(g^{\w })\neq 0$, 
since $\Phi ^{\w ,g}$ is a nonzero element of $K^{\w }[y]$.
Hence, $m_{\w }^g(\Phi )=0$ and 
$\deg \Phi (g)=\deg ^g\Phi $
by Lemma~\ref{lem:initial-equiv}(i).

(ii) Set $P_0=P(K^{\w },g^{\w })$, 
$P=P(A^{\w },g^{\w })$ and $I=I(A^{\w },g^{\w })$. 
By Lemma~\ref{lem:initial-equiv}(ii), 
we have 
$(\partial _y^{m-1}\Phi ^{\w ,g})(g^{\w })=0$ 
and $(\partial _y^{m}\Phi ^{\w ,g})(g^{\w })\neq 0$. 
Since $k$ is of characteristic zero, 
this implies that 
$\Phi ^{\w ,g}=P_0^mH$ for some $H\in K^{\w }[y]$ with $H(g^{\w })\neq 0$. 
By the assumption that 
$g^{\w }$ is algebraic over $K^{\w }$, 
it follows that 
$\deg _yP_0
=[K^{\w }(g^{\w }):K^{\w }]$. 
Thus, we get 
$\deg _y\Phi ^{\w ,g}=m_{\w }^g(\Phi )[K^{\w }(g^{\w }):K^{\w }]+\deg _yH$. 
Therefore, 
$m_{\w }^g(\Phi )$ is at most the quotient of 
$\deg _y\Phi ^{\w ,g}$ divided by $[K^{\w }(g^{\w }):K^{\w }]$. 
Assume that $I$ is a principal ideal. 
Write $\Phi ^{\w ,g}={P}^{m'}H'$, 
where $m'\in \Zn $ and 
$H'\in A^{\w }[y]\sm I$. 
Then, 
$m'$ must be at most $m$, 
since $P$ belongs to $P_0K^{\w }[y]$. 
On the other hand, 
$P$ does not belong to $P_0^2K^{\w }[y]$, 
for otherwise $\partial _yP$ 
would belong to 
$P_0K^{\w }[y]\cap A^{\w }[y]=I=PA^{\w }[y]$, 
a contradiction. 
Hence, 
$m'$ must be at least $m$, 
since $H'(g^{\w })\neq 0$. 
Thus, $m'=m$. 
This proves the latter part. 
\qed\end{pf}

Here is a generalization of 
the Shestakov-Umirbaev inequality~\cite[Theorem 3]{SU1}.

\begin{thm}\label{thm:generalSU}
Let $f_1,\ldots ,f_r$ and $g$ be nonzero 
elements of $\kx $ for $r\geq 1$ with 
$f_1,\ldots ,f_r$ algebraically independent over $k$, 
and let $A=k[f_1,\ldots ,f_r]$, $K=k(f_1,\ldots ,f_r)$ 
and $\omega =df_1\wedge \cdots \wedge df_r$. 
Let $\w \in \Gamma ^n$ 
such that $\deg _{\w }h\geq 0$ for each $h\in A\sm \zs $, 
and $M=\deg _{\w }(\omega \wedge dg)-\deg _{\w }\omega -\deg _{\w }g$. 
Then, we have the following for each $\Phi \in A[y]\sm \zs $:  

{\rm (i)} 
Assume that 
$g^{\w }$ is algebraic over $K^{\w }$, 
and let $a$ and $b$ be the quotient and residue 
of $\deg _y\Phi $ 
divided by $[K^{\w }(g^{\w }):K^{\w }]$, respectively. 
Then, 
it follows that 
\begin{equation}\label{eq:generalSU1}
\deg _{\w }\Phi (g)\geq (\deg _y\Phi )\deg _{\w }g +aM
=a\!\left( 
[K^{\w }(g^{\w }):K^{\w }]\deg _{\w }g+M\right) +b\deg _{\w }g. 
\end{equation}

{\rm (ii)} 
If $I(A^{\w },g^{\w })$ is a principal ideal and $\deg _{\w }g\geq 0$, 
then 
\begin{equation}\label{eq:generalSU2}
\deg _{\w }\Phi (g)\geq m_{\w }^g(\Phi )(\deg _{\w }^gP(A^{\w },g^{\w })+M). 
\end{equation}
\end{thm}
\begin{pf}\rm
(i) 
The equality in (\ref{eq:generalSU1}) can be checked easily. 
We only show the inequality. 
By Theorem~\ref{thm:main}, 
we get $\deg _{\w }\Phi (g)\geq \deg _{\w }^g\Phi +m_{\w }^g(\Phi )M$. 
It suffices to verify that 
$\deg _{\w }^g\Phi \geq (\deg _y\Phi )\deg _{\w }g$ 
and $m_{\w }^g(\Phi )M\geq qM$. 
Let $\phi _e\in A$ be the coefficient of $y^e$ in $\Phi $, 
where $e=\deg _y\Phi $. 
Then, 
$\deg _{\w }^g\Phi \geq \deg _{\w }(\phi _eg^e)$. 
Besides, 
$\deg _{\w }\phi _e\geq 0$ by the assumption on $\w $. 
Hence, we get 
$$
\deg _{\w }^g\Phi \geq \deg _{\w }(\phi _eg^e)
=\deg _{\w }\phi _e+e\deg _{\w }g
\geq (\deg _y\Phi )\deg _{\w }g. 
$$
On the other hand, we obtain 
$M\leq 0$ using (\ref{eq:ref0}) and (\ref{eq:ref1}). 
Moreover, 
$m_{\w }^g(\Phi )\leq a$ by Proposition~\ref{prop:m_g}(ii). 
Therefore, 
$m_{\w }^g(\Phi )M\geq aM$, 
proving the inequality in (\ref{eq:generalSU1}).

(ii) 
We note that $\deg _{\w }^g\Psi \geq 0$ 
whenever $\Psi $ is a nonzero element of 
$A[y]\cup A^{\w }[y]$. 
Actually, $\deg _{\w }^g\Psi =\deg _{\w }\psi +l\deg _{\w }g$ 
for some $\psi \in A\sm \zs $ and $l\in \Zn $, 
and $\deg _{\w }\psi \geq 0$ and $\deg _{\w }g\geq 0$ by assumption. 
First, 
assume that $g^{\w }$ is transcendental over $K^{\w }$. 
Then, $m _{\w }^g(\Phi ^{\w ,g})=0$ 
and 
$\deg _{\w }\Phi (g)=\deg _{\w }^g\Phi $ 
by Proposition~\ref{prop:m_g}(i). 
Hence, 
the right-hand side of (\ref{eq:generalSU2}) is zero, 
while $\deg _{\w }(\Phi (g))\geq 0$, 
since $\deg _{\w }^g\Phi\geq 0$ as noted. 
Therefore, 
(\ref{eq:generalSU2}) is true if 
$g^{\w }$ is transcendental over $K^{\w }$. 
Next, assume that $g^{\w }$ is algebraic over $K^{\w }$. 
By Proposition~\ref{prop:m_g}(ii), 
we get $\Phi ^{\w ,g}=P^mH$ for some $H\in A^{\w }[y]$, 
where $P=P(A^{\w },g^{\w })$ and $m=m_{\w }^g(\Phi )$. 
Since $\deg _{\w }H\geq 0$ as noted, 
we obtain 
$$
\deg _{\w }^g\Phi =\deg _{\w }^g\Phi ^{\w ,g}=
m\deg _{\w }^gP+\deg _{\w }^gH
\geq m_{\w }^g(\Phi )\deg _{\w }^gP. 
$$
With the aid of this inequality, 
(\ref{eq:generalSU2}) 
follows from Theorem~\ref{thm:main}. 
\qed\end{pf}

The following lemma is well-known. 
For the sake of completeness, 
we include a proof at the end of this section.

\begin{lem}\label{lem:fgdep}
Let $f$ and $g$ be $\w $-homogeneous elements of $\kx $ 
with $\deg _{\w }f>0$ and $\deg _{\w }g>0$ for some $\w \in \Gamma ^n$. 
If $f$ and $g$ are algebraically dependent over $k$, 
then there exist mutually prime natural numbers 
$l(f,g)$ and $l(g,f)$ as follows: 

{\rm (i)} 
$g^{l(f,g)}=\alpha f^{l(g,f)}$ for some $\alpha \in k$.

{\rm (ii)} 
$I(k[f],g)=\left( 
y^{l(f,g)}-\alpha f^{l(g,f)}\right)\!
k[f][y]$. 

{\rm (iii)} $[k(f)(g):k(f)]=l(f,g)$. 

{\rm (iv)} 
$l(f,g)=(\deg _{\w }f)\gcd (\deg _{\w }f,\deg _{\w }g)^{-1}$ 
if $\Gamma =\Z $. 
\end{lem}

The Shestakov-Umirbaev inequality~\cite[Theorem 3]{SU1} 
is obtained as a corollary to Theorem~\ref{thm:generalSU}.

\begin{cor}[Shestakov-Umirbaev]\label{cor:SU}
Assume that $f,g\in \kx \sm k$ satisfy 
$\deg _{\w }f>0$ and $\deg _{\w }g>0$ for some $\w \in \Z ^n$. 
Then, for each $\Phi \in k[f][y]\sm \zs $, 
it follows that 
\begin{align}
\begin{split}\label{eq:corSU}
\deg _{\w }\Phi (g)
\geq 
a(\lcm (\deg _{\w }f,\deg _{\w }g)
+M)
+b\deg _{\w }g
\end{split}
\end{align}
where $M=\deg _{\w }(df\wedge dg)-\deg _{\w }f-\deg _{\w }g$, 
and $a$ and $b$ are the quotient and residue of $\deg _y\Phi $ 
divided by $(\deg _{\w }f)\gcd (\deg _{\w }f,\deg _{\w }g)^{-1}$, 
respectively. 
\end{cor}
\begin{pf}\rm
We remark that 
$k[f]^{\w }=k[f^{\w }]$, 
and $\deg _{\w }h\geq 0$ for each $h\in k[f]\sm \zs $. 
In fact, 
if $h=\sum _{i=0}^ec_if^i$, where 
$c_0,\ldots ,c_e\in k$ with $c_e\neq 0$ for $e\geq 0$, 
then $\deg _{\w }h=e\deg _{\w }f\geq 0$ 
and $h^{\w }=c_e(f^{\w })^e$, 
since $\deg _{\w }f>0$ by assumption. 
Consequently, we have 
$k(f)^{\w }=k(f^{\w })$. 
First, 
assume that $f^{\w }$ and $g^{\w }$ are algebraically dependent over $k$, 
and put $N=[k(f^{\w })(g^{\w }):k(f^{\w })]$. 
Then, 
Theorem~\ref{thm:generalSU}(i) gives that 
\begin{equation}\label{eq:corSU'}
\deg _{\w }\Phi (g)\geq a'(N\deg _{\w }g+M)+b'\deg _{\w }g, 
\end{equation}
where $a'$ and $b'$ are the quotient and residue of $\deg _y\Phi $ 
divided by $N$, 
respectively. 
By Lemma~\ref{lem:fgdep}, 
we have 
\begin{equation*}
N=\frac{\deg _{\w }f^{\w }}{\gcd(\deg _{\w }f^{\w },\deg _{\w }g^{\w })}
=\frac{\deg _{\w }f}{\gcd(\deg _{\w }f,\deg _{\w }g)}
=\frac{\lcm (\deg _{\w }f,\deg _{\w }g)}{\deg _{\w }g}. 
\end{equation*}
This implies that the right-hand side of (\ref{eq:corSU'}) 
is equal to that of (\ref{eq:corSU}). 
Therefore, (\ref{eq:corSU}) is true. 
If $f^{\w }$ and $g^{\w }$ are algebraically independent over $k$, 
then $\deg _{\w }\Phi (g)=\deg _{\w }^g\Phi $ 
by Proposition~\ref{prop:m_g}(i). 
As in the proof of Theorem~\ref{thm:generalSU}, 
we get 
$\deg _{\w }^g\Phi \geq (\deg _y\Phi )\deg _{\w }g$. 
On the other hand, 
the right-hand side of (\ref{eq:corSU}) 
is equal to 
$(\deg _y\Phi )\deg _{\w }g+aM$, 
and also $M\leq 0$. 
This proves (\ref{eq:corSU}). 
\qed\end{pf}

In the original statement of~\cite[Theorem 3]{SU1}, 
the ``Poisson bracket" $[f,g]$ is used instead of $df\wedge dg$. 
The degrees of $[f,g]$ and $df\wedge dg$ are defined in the same way.

To conclude this section, 
we prove Lemma~\ref{lem:fgdep}. 
The assertions (ii), (iii) and (iv) 
easily follows from the assertion (i). 
We only show that there exist mutually prime natural numbers 
$l$ and $m$ such that $f^{-m}g^l$ belongs to $k$. 
Without loss of generality, 
we may assume that $k$ is algebraically closed. 
In fact, 
$f^{-m}g^l$ necessarily belongs to $k$ 
if $f^{-m}g^l$ is algebraic over $k$, 
since the field of fractions of $\kx $ is a regular extension of $k$.

By the assumption that $f$ and $g$ 
are algebraically dependent over $k$, 
we may find a nontrivial algebraic relation 
$\sum _{i,j}\beta _{i,j}f^ig^j=0$, 
where $\beta _{i,j}\in k$ for each $i,j\in \Zn $. 
Let $J$ be the set of $(i,j)\in (\Zn )^2$ 
such that $\beta _{i,j}\neq 0$, 
and $(i_0,j_0)$ and $(i_1,j_1)$ the elements of $J$ 
such that $i_0\leq i\leq i_1$ 
for each $i\in \Zn $ with $(i,j)\in J$ for some $j$. 
Since $f$ and $g$ are $\w $-homogeneous, 
we may assume that $i\deg _{\w }f+j\deg _{\w }g$ 
are the same for any $(i,j)\in J$. 
Then, 
$(i_1-i_0)\deg _{\w }g=(j_0-j_1)\deg _{\w }f$. 
We note that $i_1-i_0$ must be positive, 
for otherwise $J=\{ (i_0,j_0)\} $, 
and then 
$0=\sum _{(i,j)\in J}\beta _{i,j}f^ig^j
=\beta _{i_0,j_0}f^{i_0}g^{j_0}\neq 0$, 
a contradiction. 
Since 
$\deg _{\w }f>0$ and $\deg _{\w }g>0$ by assumption, 
we get $j_0-j_1>0$. 
Set $l'=i_1-i_0$, $m'=j_0-j_1$ 
and $l=l'/e$, $m=m'/e$, 
where $e=\gcd (l',m')$. 
Then, 
$J$ is contained in the set of 
$(i_0,j_0)+p(l,-m)$ for $p=0,\ldots ,e$. 
By putting $\beta _p'=\beta _{i_0+lp,j_0-mp}$ for each $p$, 
we get
$$
0=\sum _{(i,j)\in J}\beta _{i,j}f^ig^j
=f^{j_0}g^{i_0}\sum _{p=0}^e\beta _p'(f^{-m}g^{l})^p 
=\beta _e'f^{j_0}g^{i_0}\prod _{p=1}^e(f^{-m}g^{l}-\alpha _p), 
$$
where 
$\alpha _1,\ldots ,\alpha _e\in k$ 
are the solutions of the algebraic equation 
$\sum _{p=0}^e\beta _p'y^p=0$ in $y$. 
Thus, 
$f^{-m}g^{l}=\alpha _p$ for some $p$. 
Therefore, 
$f^{-m}g^{l}$ is contained in $k$. 
This completes the proof of Lemma~\ref{lem:fgdep}.

\section{A characterization of polynomial automorphisms}
\setcounter{equation}{0}
\label{sect:auto}

As an application of our result, 
we study features of elements of $\Aut _k\kx $. 
Namely, 
we give a characterization of $n$-tuples 
${\bf f}=(f_1,\ldots ,f_n)$ 
of elements of $\kx $ 
such that $k[f_1,\ldots ,f_n]=\kx $.

First, 
we recall a basic fact about initial algebras.

\begin{lem}\label{lem:initial}
Let $g_1,\ldots ,g_r$ be elements of $\kx $ 
for $r\geq 0$. 
If $g_1^{\w },\ldots ,g_r^{\w }$ 
are algebraically independent over $k$ 
for $\w \in \Gamma ^n$, 
then $k[g_1,\ldots ,g_r]^{\w }
=k[g_1^{\w },\ldots ,g_r^{\w }]$. 
\end{lem}
\begin{pf}\rm
Clearly, 
$k[g_1,\ldots ,g_r]^{\w }$ contains $k[g_1^{\w },\ldots ,g_r^{\w }]$. 
We show the reverse inclusion by induction on $r$. 
The assertion is obvious if $r=0$. 
Assume that $r\geq 1$. 
It suffices to verify that 
$h^{\w }$ belongs to $k[g_1^{\w },\ldots ,g_r^{\w }]$ 
for each $h\in k[g_1,\ldots ,g_r]\sm \zs $. 
Take $H\in A[y]$ 
such that $h=H(g_r)$, 
where $A=k[g_1,\ldots ,g_{r-1}]$. 
By induction assumption, we have 
$A^{\w }=k[g_1^{\w },\ldots ,g_{r-1}^{\w }]$. 
Besides, 
$H^{\w ,g_r}$ belongs to $A^{\w }[y]\sm \zs $. 
Hence, 
$H^{\w ,g_r}(g_r^{\w })$ 
is contained in $k[g_1^{\w },\ldots ,g_r^{\w }]$. 
Moreover, $H^{\w ,g_r}(g_r^{\w })$ is not zero, 
since $g_1^{\w },\ldots ,g_r^{\w }$ 
are algebraically independent over $k$ by assumption. 
Hence, 
$H(g_r)^{\w }=H^{\w ,g_r}(g_r^{\w })$ 
by Lemma~\ref{lem:initial-equiv}(i). 
Since $h=H(g_r)$, we get 
$h^{\w }=H(g_r)^{\w }$. 
Thus, 
$h^{\w }$ belongs to $k[g_1^{\w },\ldots ,g_r^{\w }]$. 
Therefore,  
$k[g_1,\ldots ,g_r]^{\w }$ is contained in 
$k[g_1^{\w },\ldots ,g_r^{\w }]$. 
\qed\end{pf}

The following proposition 
is an immediate consequence of Lemma~\ref{lem:initial}.

\begin{prop}\label{prop:linear}
Let $f_1,\ldots ,f_n$ be elements of $\kx $ 
such that $k[f_1,\ldots ,f_n]=\kx $. 
Then, $f_1^{\w },\ldots ,f^{\w }_n$ 
are algebraically independent over $k$ 
if and only if $k[f^{\w }_1,\ldots ,f^{\w }_n]=\kx $ 
for $\w \in \Gamma ^n$. 
\end{prop}
\begin{pf}\rm
The ``if" part is clear, 
for $\kx $ has transcendence degree $n$ over $k$. 
Assume that $f^{\w }_1,\ldots ,f^{\w }_n$ 
are algebraically independent over $k$. 
Then, 
$k[f^{\w }_1,\ldots ,f^{\w }_n]=k[f_1,\ldots ,f_n]^{\w }$ 
by Lemma~\ref{lem:initial}. 
Since $k[f_1,\ldots ,f_n]=\kx $, 
we have $k[f_1,\ldots ,f_n]^{\w }=\kx ^{\w }=\kx $. 
Thus, 
$k[f^{\w }_1,\ldots ,f^{\w }_n]=\kx $. 
This proves the ``only if" part. 
\qed\end{pf}

Next, we consider the case where 
$k(f_1^{\w },\ldots ,f_n^{\w })$ 
has transcendence degree $n-1$ over $k$ 
for some $\w \in \Gamma ^n$. 
We define an element $\Delta _{\bf f}^{\w }$ of $\Gamma $ as follows: 
Let $\lambda _{\bf f}^{\w }:\kx \to \kx $ 
be the homomorphism defined by $\lambda (x_i)=f^{\w }_i$ 
for $i=1,\ldots ,n$. 
Then, 
$\ker \lambda _{\bf f}^{\w }$ 
is a prime ideal of $\kx $ of hight one. 
Since 
$\kx $ is a unique factorization domain, 
there exists $Q\in \kx \sm \zs $ 
such that $\ker \lambda _{\bf f}^{\w }=Q\kx $. 
We define $\Delta _{\bf f}^{\w }$ to be the 
$\w _{\bf f}$-degree of $Q$, 
where 
$$
\w _{\bf f}=(\deg _{\w }f_1,\ldots ,\deg _{\w }f_n). 
$$ 
Note that $\Delta _{\bf f}^{\w }$ 
is uniquely determined by ${\bf f}$ and $\w $, 
since $Q$ is unique up to multiplication by elements in $k\sm \zs $.

Here is the main theorem of this section.

\begin{thm}\label{thm:general-Jung}
Let $f_1,\ldots ,f_n$ be elements of $\kx $ 
such that $k[f_1,\ldots ,f_n]=\kx $, 
and $\w =(w_1,\ldots ,w_n)$ an element of $(\Gamma _{\geq 0})^n$. 
If $k(f^{\w }_1,\ldots ,f^{\w }_n)$ has transcendence degree $n-1$ over $k$, 
then 
\begin{equation}\label{eq:general-jung}
\sum _{i=1}^n\deg _{\w }f_i
\geq \Delta _{\bf f}^{\w }+\sum _{i=1}^nw_i
-\max \{ w_i\mid i=1,\ldots ,n\} ,
\end{equation}
where ${\bf f}=(f_1,\ldots ,f_n)$. 
\end{thm}
\begin{pf}\rm
Since $k(f^{\w }_1,\ldots ,f^{\w }_n)$ 
has transcendence degree $n-1$ over $k$, 
we may find $l$ such that $x_l$ 
is not contained in $k[f^{\w }_1,\ldots ,f^{\w }_n]$. 
Moreover, 
we may assume that $f^{\w }_1,\ldots ,f^{\w }_{n-1}$ 
are algebraically independent over $k$ 
by changing the indices of $f_1,\ldots ,f_n$ if necessary. 
Set $A=k[f_1,\ldots ,f_{n-1}]$ and $g=f_n$. 
Then, 
there exists $\Phi \in A[y]$ such that $\Phi (g)=x_l$, 
since $A[g]=\kx $ by assumption. 
Furthermore, 
$A^{\w }=k[f_1^{\w },\ldots ,f_{n-1}^{\w }]$ 
by Lemma~\ref{lem:initial}, 
and so $A^{\w }$ is a polynomial ring over $k$. 
Accordingly, 
$I(A^{\w },g^{\w })$ is a principal ideal of $A^{\w }[y]$. 
Besides, 
$\deg _{\w }h\geq 0$ holds for each $h\in \kx \sm \zs $, 
since $w_i\geq 0$ for $i=1,\ldots ,n$ by assumption. 
Then, we can easily check that 
$f_1,\ldots ,f_{n-1}$, $g$ and $\w $ satisfy 
the assumptions of Theorem~\ref{thm:generalSU}(ii). 
Therefore, we obtain 
\begin{equation}\label{eq:generaljung2}
\deg _{\w }\Phi (g)\geq m_{\w }^g(\Phi )
(\deg _{\w }^gP+M), 
\end{equation}
where $P=P(A^{\w },g^{\w })$, 
$\omega =df_1\wedge \cdots \wedge df_{n-1}$ and 
$M=\deg _{\w }(\omega \wedge dg)-\deg _{\w }\omega -\deg _{\w }g$. 
We show that 
\begin{equation}\label{eq:Mgeq}
M\geq \sum _{i=1}^nw_i-\sum _{i=1}^n\deg _{\w }f_i. 
\end{equation}
Note that 
$\omega \wedge dg=df_1\wedge \cdots \wedge df_n=
\alpha dx_1\wedge \cdots \wedge dx_n$, 
where $\alpha $ is the determinant of the $n$ by $n$ matrix 
$(\partial f_i/\partial x_j)_{i,j}$. 
The assumption $k[f_1,\ldots ,f_n]=\kx $ 
implies that $\alpha $ belongs to $k\sm \zs $. 
Hence, 
we have 
\begin{equation}\label{eq:degM1}
\deg _{\w }(\omega \wedge dg)
=\deg _{\w }(\alpha dx_1\wedge \cdots \wedge dx_n)
=\deg _{\w }\alpha +\sum _{i=1}^nw_i
=\sum _{i=1}^nw_i. 
\end{equation}
On the other hand, 
we get 
\begin{equation}\label{eq:degM2}
\deg _{\w }\omega 
=\deg _{\w }(df_1\wedge \cdots \wedge df_{n-1})
\leq \sum _{i=1}^{n-1}\deg _{\w }df_i
=\sum _{i=1}^{n-1}\deg _{\w }f_i 
\end{equation}
by using (\ref{eq:ref0}) and (\ref{eq:ref1}). 
Since $g=f_n$, 
the inequality (\ref{eq:Mgeq}) follows from 
(\ref{eq:degM1}) and (\ref{eq:degM2}).

To complete the proof, 
it remains only to show that 
$m_{\w }^g(\Phi )\geq 1$ and $\deg _{\w }^gP=\Delta _{\bf f}^{\w }$. 
Actually, assuming this, 
we can easily deduce (\ref{eq:general-jung}) 
from the inequalities (\ref{eq:generaljung2}), (\ref{eq:Mgeq}) 
and 
$$
\max \{ w_i\mid i=1,\ldots ,n\} \geq w_l=\deg _{\w }x_l=
\deg _{\w }\Phi (g). 
$$
First, 
suppose to the contrary that $m_{\w }^g(\Phi )=0$. 
Then, 
$\Phi ^{\w ,g}(g^{\w })=\Phi (g)^{\w }=x_l^{\w }=x_l$ 
by Lemma~\ref{lem:initial-equiv}. 
Recall that $x_l$ does not belong to 
$k[f_1^{\w },\ldots ,f_n^{\w }]$, 
while 
$k[f_1^{\w },\ldots ,f_n^{\w }]=A^{\w }[g^{\w }]$. 
Since $\Phi ^{\w ,g}$ is in $A^{\w }[y]$, 
it follows that 
$\Phi ^{\w ,g}(g^{\w })$ belongs to $A^{\w }[g^{\w }]$. 
This is a contradiction. 
Thus, we get $m_{\w }^g(\Phi )\geq 1$. 
Next, take $Q\in \kx $ so that 
$\ker \lambda _{\bf f}^{\w }=Q\kx $. 
Let $\iota :\kx \to A^{\w }[y]$ be 
the homomorphism defined by 
$\iota (x_i)=f_i^{\w }$ 
for $i=1,\ldots ,n-1$ and $\iota (x_n)=y$. 
Then, $\iota $ is an isomorphism, 
since we are assuming that $f_1^{\w },\ldots ,f_{n-1}^{\w }$ 
are algebraically independent over $k$. 
This assumption implies further that 
the $\w _{\bf f}$-degree of $Q$ is equal to 
the $(\w ,\deg _{\w }^g)$-degree of $\iota (Q)$. 
It is equal to $\deg _{\w }^g\iota (Q)$ 
as mentioned at the beginning of Section 3. 
Thus, 
we get $\Delta _{\bf f}^{\w }=\deg _{\w }^g\iota (Q)$. 
By definition, 
$\lambda _{\bf f}^{\w }$ is equal to the composite of 
$\iota $ and the substitution map 
$A^{\w }[y]\ni \psi \mapsto \psi (g^{\w })\in \kx $. 
Hence, we have 
$$
\iota (Q\kx )=\iota (\ker \lambda _{\bf f}^{\w })
=I(A^{\w },g^{\w })=PA^{\bf w}[y]. 
$$
Since $\iota $ is an isomorphism, 
$\iota (Q)=\alpha P$ for some $\alpha \in k\sm \zs $. 
Thus, 
$\deg _{\w }^g\iota (Q)=\deg _{\w }^gP$. 
Therefore, 
we obtain $\Delta _{\bf f}^{\w }=\deg _{\w }^gP$. 
\qed\end{pf}

Theorem~\ref{thm:general-Jung}
is considered as a generalization of Proposition~\ref{prop:Jung}. 
In fact, 
we have the following corollary in case $n=2$.

\begin{cor}\label{cor:Jung}
Assume that $f_1,f_2\in k[x_1,x_2]$ satisfy 
$k[f_1,f_2]=k[x_1,x_2]$. 
If $f_1^{\w }$ and $f_2^{\w }$ are algebraically dependent over $k$ 
for $\w \in (\Zn )^2$, 
then $\deg _{\w }f_1$ and $\deg _{\w }f_2$ are positive integers 
which satisfy 
\begin{equation}\label{eq:cor:Jung}
\deg _{\w }f_1+\deg _{\w }f_2\geq 
\lcm (\deg _{\w }f_1,\deg _{\w }f_2)+\min \{ w_1,w_2\} , 
\end{equation}
where $\w =(w_1,w_2)$. 
In particular, 
$\deg _{\w }f_1|\deg _{\w }f_2$ or $\deg _{\w }f_2|\deg _{\w }f_1$. 
\end{cor}
\begin{pf}\rm
Since $w_i\geq 0$ for $i=1,2$ by assumption, 
$\deg _{\w }f_i\geq 0$ for $i=1,2$. 
We show that 
$\deg _{\w }f_i\neq 0$ for $i=1,2$ by contradiction. 
Suppose the contrary, say $\deg _{\w }f_1=0$. 
Then, $w_i=0$ for some $i\in \{ 1,2\} $, 
since $f_1$ cannot be contained in $k$. 
We claim that $\w \neq 0$, 
for otherwise $f_i^{\w }=f_i$ for $i=1,2$. 
This is impossible, 
because $k[f_1,f_2]=k[x_1,x_2]$, 
whereas $f_1^{\w }$ and $f_2^{\w }$ 
are algebraically dependent over $k$. 
Hence, 
we have $w_j>0$ for $j\in \{ 1,2\} \sm \{ i\} $. 
Since we suppose that $\deg _{\w }f_1=0$, 
this implies that $f_1$ belongs to $k[x_i]$, 
and besides $f_1^{\w }=f_1$. 
Then, $f_2^{\w }$ also belongs to $k[x_i]$, 
since $f_1^{\w }$ and $f_2^{\w }$ are algebraically dependent over $k$. 
Consequently, 
$f_2$ belongs to $k[x_i]$ 
due to the conditions $w_i=0$ and $w_j>0$. 
Thus, 
$k[f_1,f_2]$ is contained in $k[x_i]$, 
a contradiction. 
Therefore, 
$\deg _{\w }f_i\neq 0$ for $i=1,2$.

Put $P=P(k[f_1^{\w }], f_2^{\w })$ and ${\bf f}=(f_1,f_2)$. 
As in the proof of Theorem~\ref{thm:general-Jung}, 
we have 
$\Delta _{\bf f}^{\w }=\deg _{\w }^{f_2}P$. 
By Lemma~\ref{lem:fgdep}, 
we may write 
$P=\beta\!\left( 
y^{l(f_1,f_2)}-\alpha (f_1^{\w })^{l(f_2,f_1)}\right) $, 
where $\alpha ,\beta \in k\sm \zs $. 
Then, we have 
$\deg _{\w }^{f_2}P=\lcm (\deg _{\w }f_1,\deg _{\w }f_2)$. 
Thus, 
$\Delta _{\bf f}^{\w }=\lcm (\deg _{\w }f_1,\deg _{\w }f_2)$. 
By Theorem~\ref{thm:general-Jung}, we obtain 
\begin{equation*}
\deg _{\w }f_1+\deg _{\w }f_2\geq 
\Delta _{\bf f}^{\w }+w_1+w_2-\max \{ w_1,w_2\} 
=\lcm (\deg _{\w }f_1,\deg _{\w }f_2)+\min \{ w_1,w_2\} 
\end{equation*}
The last statement is a consequence of 
the first statement, 
since 
$a+b\geq \lcm (a,b)$ implies 
$a|b$ or $b|a$ for each $a,b\in \N $. 
\qed
\end{pf}

\section{A lemma of Shestakov-Umirbaev}
\setcounter{equation}{0}

For $f_1,f_2,f_3\in \kx \sm k$, we put 
$$
m_1=\deg f_1+\deg _{\w }(df_2\wedge df_3),\ 
m_2=\deg f_2+\deg _{\w }(df_3\wedge df_1),\ 
m_3=\deg f_3+\deg _{\w }(df_1\wedge df_2), 
$$
where $\w =(1,\ldots ,1)$. 
Shestakov-Umirbaev~\cite[Lemma 5]{SU1} 
proved the following lemma. 

\begin{lem}[Shestakov-Umirbaev]\label{lem:jacobi}
In the notation above, 
$m_1\leq \max \{ m_2,m_3\} $. 
If $m_2\neq m_3$, then $m_1=\max \{ m_2,m_3\} $. 
\end{lem}
This lemma also plays an important role in~\cite{SU2} 
to solve the Nagata conjecture. 
We note that the statement of Lemma~\ref{lem:jacobi} 
is equivalent to the following statement: 

($\dag $) There exist $1\leq i_1<i_2\leq 3$ 
such that $m_{i_1}=m_{i_2}\geq m_i$ for $i=1,2,3$. \\
To conclude this paper, 
we give a generalization of the lemma of Shestakov-Umirbaev. 

\begin{thm}\label{thm:jacobi}
Let $\eta _1,\ldots ,\eta _{l}$ be elements of $\Omega _{\kx /k}$ 
for $l\geq 2$. 
Then, there exist $1\leq i_1<i_2\leq l$ 
such that 
$$
\deg _{\w }\eta _{i_1}+\deg _{\w }\tilde{\eta }_{i_1}
=\deg _{\w }\eta _{i_2}+\deg _{\w }\tilde{\eta }_{i_2}
\geq \deg _{\w }\eta _{i}+\deg _{\w }\tilde{\eta }_{i}
$$
for $i=1,\ldots ,l$, 
where $\tilde{\eta }_i
=\eta _1\wedge \cdots \wedge \eta _{i-1}\wedge 
\eta _{i+1}\wedge \cdots \wedge \eta _{l}$ for each $i$. 
\end{thm}

The statement ($\dag $) 
is obtained as a special case of Theorem~\ref{thm:jacobi} 
where $l=3$ and $\eta _i=df_i$ for $i=1,2,3$, 
since
$$
\deg _{\w }\eta _i+\deg _{\w }\tilde{\eta }_i
=\deg _{\w }df_i+\deg _{\w }(df_j\wedge df_k)
=\deg _{\w }f_i+\deg _{\w }(df_j\wedge df_k)=m_i$$ 
for each distinct integers $1\leq i,j,k\leq 3$ in this case.

Let us prove Theorem~\ref{thm:jacobi} by contradiction. 
Suppose to the contrary that 
there exists $i_0$ such that 
$\deg _{\w }\eta _{i_0}+\deg _{\w }\tilde{\eta }_{i_0}
>\deg _{\w }\eta _{i}+\deg _{\w }\tilde{\eta }_{i}$ 
for each $i\neq i_0$. 
Write $\eta _i=\sum _{j=1}^nf_{i,j}x_j^{-1}dx_j$ for each $i$, 
where $f_{i,j}\in x_j\kx $ for each $j$. 
Set $d\x _I=dx_{i_1}\wedge \cdots \wedge dx_{i_{l-1}}$ 
and $\x _I=x_{i_1}\cdots x_{i_{l-1}}$ 
for each $i_1,\ldots ,i_{l-1}$, where $I=(i_1,\ldots ,i_{l-1})$. 
Then, 
we may write 
$\tilde{\eta }_i=\sum _J\tilde{f}_{i,J}(\x _J)^{-1}d\x _J$, 
where the sum is taken over $J=(j_1,\ldots ,j_{l-1})$ 
with $1\leq j_1<\cdots <j_{l-1}\leq n$, 
and $\tilde{f}_{i,J}\in \x _J\kx $ for each $J$. 
By the definition (\ref{eq:def:deg}) of the $\w $-degree, 
there exist $j_0$ and $J_0$ such that 
$\deg _{\w }\eta _{i_0}=\deg _{\w }f_{i_0,j_0}$ and 
$\deg _{\w }\tilde{\eta }_{i_0}=\deg _{\w }\tilde{f}_{i_0,J_0}$. 
By the choice of $i_0$, 
it follows that 
$\deg _{\w }(f_{i,j}\tilde{f}_{i,J})<
\deg _{\w }(f_{i_0,j_0}\tilde{f}_{i_0,J_0})$ 
for each $j$ and $J$ if $i\neq i_0$. 
In particular, 
$f_{i_0,j_0}\neq 0$ and $\tilde{f}_{i_0,J_0}\neq 0$. 
By changing the indices of 
$\eta _1,\ldots ,\eta _l$ and 
$x_1,\ldots ,x_n$ if necessary, 
we may assume that $i_0\neq 1$ and $J_0=(1,\ldots ,l-1)$. 
Note that the $(i,l)$-cofactor of 
the $l$ by $l$ matrix 
$$
M=\left(
\begin{array}{@{\!}cccc@{\!}}
f_{1,1}&\cdots & f_{1,l-1}& f_{1,j_0}\\
f_{2,1}&\cdots & f_{2,l-1}& f_{2,j_0}\\
 & \hdotsfor{2} & \\
f_{l,1}&\cdots & f_{l,l-1}& f_{l,j_0}
\end{array}
\right) 
$$
is equal to $(-1)^{l+i}\tilde{f}_{i,J_0}$ for $i=1,\ldots ,l$. 
Hence, 
$\det M=\sum _{i=1}^{l}(-1)^{i}f_{i,j_0}\tilde{f}_{i,J_0}$. 
Since 
$\deg _{\w }(f_{i,j_0}\tilde{f}_{i,J_0})
<\deg _{\w }(f_{i_0,j_0}\tilde{f}_{i_0,J_0})$ 
if $i\neq i_0$, 
we get 
$\deg _{\w }(\det M)=\deg _{\w }(f_{i_0,j_0}\tilde{f}_{i_0,J_0})$. 
On the other hand, 
the $(1,u)$-cofactor of $M$ is equal to 
$(-1)^{u}\tilde{f}_{1,J_u}$ for $u=1,\ldots ,l$, 
where $J_u=(1,\ldots ,u-1,u+1,\ldots ,l-1,j_0)$ 
for $1\leq u<l$ and $J_l=J_0$. 
Hence, 
$\det M=\sum _{u=1}^{l}(-1)^{u}f_{1,u}\tilde{f}_{1,J_u}$. 
Since we assume that $i_0\neq 1$, 
it follows that 
$\deg _{\w }(f_{1,u}\tilde{f}_{1,J_u})
<\deg _{\w }(f_{i_0,j_0}\tilde{f}_{i_0,J_0})$ 
for each $u$. 
Thus, 
$\deg _{\w }(\det M)<\deg _{\w }(f_{i_0,j_0}\tilde{f}_{i_0,J_0})$, 
and we are led to a contradiction. 
This completes the proof of Theorem~\ref{thm:jacobi}.

\begin{flushleft}
Department of Mathematics and Information Sciences\\ 
Tokyo Metropolitan University\\
1-1  Minami-Ohsawa, Hachioji \\
Tokyo 192-0397, Japan\\
E-mail: {\tt kuroda@tmu.ac.jp}
\end{flushleft}

\end{document}